\documentclass[a4paper,12pt]{article}

 \usepackage{latexsym,amssymb,amsmath}
 
 \newtheorem{thm}{Theorem}
 
 \newtheorem{prop}{Proposition}

\begin{document}

\title{On a sequence of monogenic polynomials satisfying the Appell condition whose first term is a non-constant function}

\author{Dixan Pe\~na Pe\~na\\\normalsize{Department of Mathematics, University of Aveiro}\\\normalsize{3810-193 Aveiro, Portugal}\\\normalsize{e-mail: dixanpena@ua.pt; dixanpena@gmail.com}}

\date{} 

\maketitle

\begin{abstract}
\noindent In this paper we aim at constructing a sequence $\{\mathsf{M}_n^k(x)\}_{n\ge0}$ of $\mathbb R_{0,m}$-valued polynomials which are monogenic in $\mathbb R^{m+1}$ satisfying the Appell condition (i.e. the hypercomplex derivative of each polynomial in the sequence equals, up to a multiplicative constant, its preceding term) but whose first term $\mathsf{M}_0^k(x)=\mathbf{P}_k(\underline x)$ is a $\mathbb R_{0,m}$-valued homogeneous monogenic polynomial in $\mathbb R^m$ of degree $k$ and not a constant like in the classical case. The connection of this sequence with the so-called Fueter's theorem will also be discussed.\vspace{0.2cm}\vspace{0.1cm}
\end{abstract}

\section{Preliminaries}

Let $\mathbb{R}_{0,m}$ be the real Clifford algebra generated by the orthonormal basis $\{e_1,\ldots,e_m\}$ of the Euclidean space $\mathbb R^m$ (see \cite{Cl}). The multiplication in $\mathbb{R}_{0,m}$ is associative and is determined by the relations:
\begin{alignat*}{2}
e_j^2=-1&,&\qquad &j=1,\dots,m,\\
e_je_k+e_ke_j=0&,&\qquad &1\le j\neq k\le m.
\end{alignat*} 
A basis for $\mathbb{R}_{0,m}$ is given by $\{e_A:\,A\subset\{1,\dots,m\}\}$ where $e_A=e_{j_1}\dots e_{j_k}$ with $A=\{j_1,\dots,j_k\}$ and $1\le j_1<\cdots<j_k\le m$. For the empty set $\emptyset$, we put $e_{\emptyset}=1$, the latter being the identity element. A general element $a\in\mathbb R_{0,m}$ may thus be written as 
\[a=\sum_Aa_Ae_A,\quad a_A\in\mathbb R,\]
and its conjugate $\overline a$ is given by
\[\overline a=\sum_Aa_A\overline e_A,\quad\overline e_A=\overline e_{j_k}\dots\overline e_{j_1},\quad\overline e_j=-e_j,\quad j=1,\dots,m.\]
Observe that $\mathbb R^{m+1}$ may be naturally embedded in the Clifford algebra $\mathbb R_{0,m}$ by associating to any element $x=(x_0,x_1,\ldots,x_m)\in\mathbb R^{m+1}$ the \emph{paravector} $x_0+\underline x=x_0+\sum_{j=1}^mx_je_j$.

Let us recall that an $\mathbb R_{0,m}$-valued function $f$ defined and continuously differentiable in an open set $\Omega$ of $\mathbb R^{m+1}$, is said to be (left) monogenic in $\Omega$ if and only if $\partial_xf(x)=0$ in $\Omega$, where
\[\partial_x=\partial_{x_0}+\partial_{\underline x}\]
is the generalized Cauchy-Riemann operator in $\mathbb R^{m+1}$ and 
\[\partial_{\underline x}=\sum_{j=1}^me_j\partial_{x_j}\]
is the Dirac operator in $\mathbb R^m$. Similarly, the same name is used for $\mathbb R_{0,m}$-valued functions defined in open subsets of $\mathbb R^m$ which are null-solutions of the Dirac operator $\partial_{\underline x}$. Note that $\partial_x$ factorizes the Laplacian, i.e. 
\[\Delta_x=\sum_{j=0}^m\partial_{x_j}^2=\partial_x\overline\partial_x=\overline\partial_x\partial_x,\]
and therefore every monogenic function is also harmonic. The monogenic functions are a fundamental object of study in Clifford analysis; and may be considered as a natural generalization to higher dimensions of the holomorphic functions of one complex variable (see e.g. \cite{BDS,DSS,GuSp}). 

The hypercomplex derivative of a monogenic function $f$ is defined as $\frac{1}{2}\,\overline\partial_xf$ (see \cite{GuM,M}). As a monogenic function $f$ clearly satisfies
\[\partial_{x_0}f=-\partial_{\underline x}f,\]
it easily follows that
\[\frac{1}{2}\,\overline\partial_xf=\partial_{x_0}f=-\partial_{\underline x}f.\]
Let us recall that a sequence of polynomials in which the index of each polynomial equals its degree is called a polynomial sequence. A polynomial sequence $\{p_n(z)\}_{n\ge0}$ is said to be an Appell sequence if it satisfies
\[p_n^\prime(z)=np_{n-1}(z),\quad n\ge1,\]
and $p_0(z)$ is a non-zero constant (see \cite{A}). Apart from the trivial example $\{z^n\}_{n\ge0}$, there are important sequences in Mathematics which are Appell sequences for example the Bernoulli polynomials, the Hermite polynomials, and the Euler polynomials. Recently, this concept has been generalized to the Clifford analysis setting in \cite{FCM,FM,MF} (see also \cite{BG,BGLS,CM,NGu,Lav}) as follows. A sequence $\{\mathsf{P}_n(x)\}_{n\ge0}$ of $\mathbb R_{0,m}$-valued polynomials forms an Appell sequence if
\begin{itemize}
\item [{\rm (i)}] $\{\mathsf{P}_n(x)\}_{n\ge0}$ is a polynomial sequence;
\item [{\rm (ii)}] each $\mathsf{P}_n(x)$ is monogenic in $\mathbb R^{m+1}$, i.e. $\partial_x\mathsf{P}_n(x)=0$ in $\mathbb R^{m+1}$;
\item [{\rm (iii)}] $\frac{1}{2}\,\overline\partial_x\mathsf{P}_n(x)=n\mathsf{P}_{n-1}(x)$, $n\ge1$.
\end{itemize}
In \cite{FCM,FM,MF} an important example of an Appell sequence of monogenic polynomials $\{\mathsf{P}_n^m(x)\}_{n\ge0}$ with $\mathsf{P}_0^m(x)=1$ was constructed, in which each term was of the form
\begin{equation}\label{ASMal}
\mathsf{P}_n^m(x)=\sum_{j=0}^n\binom{n}{j}C_n(j)x_0^j\underline x^{n-j}=\sum_{j=0}^nT_j^n(x_0+\underline x)^{n-j}(x_0-\underline x)^j,\quad n\ge0,
\end{equation}
for suitable real numbers $C_n(j)$, $T_j^n$ (see \cite{FCM,FM,MF}). The importance of these monogenic polynomials lies in the fact that they may be seen as the higher dimensional counterpart of the complex monomials $z^n$. At this point we must remark that that natural powers $(x_0+\underline x)^n$ of the paravector variable are not monogenic for $m\ge2$.  

Note that the requirement of $\{\mathsf{P}_n(x)\}_{n\ge0}$ being a polynomial sequence implies that the first term $\mathsf{P}_0(x)$ must be a constant. It is natural to ask whether we can consider sequences of monogenic polynomials satisfying the Appell condition (iii) but in which the first term is not a constant. In other words, we wish to drop condition (i) and keep conditions (ii)-(iii). The main goal of this note is to construct a sequence $\{\mathsf{M}_n^k(x)\}_{n\ge0}$ of $\mathbb R_{0,m}$-valued polynomials which are monogenic in $\mathbb R^{m+1}$ satisfying
\begin{equation}\label{Acond}
\frac{1}{2}\,\overline\partial_x\mathsf{M}_n^k(x)=n\mathsf{M}_{n-1}^k(x),\quad n\ge1,
\end{equation}
with $\mathsf{M}_0^k(x)=\mathbf{P}_k(\underline x)$ being a given but arbitrary $\mathbb R_{0,m}$-valued homogeneous polynomial of degree $k$ which is monogenic in $\mathbb R^m$.

We believe that the sequence $\{\mathsf{M}_n^k(x)\}_{n\ge0}$ provides a remarkable generalization of the well-known Appell sequences and may be of interest for those working in the subject. As far as we know, the methods we use to construct this sequence have never been employed in this framework. They allow us to arrive at our main result in a simple and elegant way. 

In a forthcoming paper, we plan to consider even more general initial functions $\mathsf{M}_0^k(x)$ such as $\mathbb R_{0,m}$-valued homogeneous polynomials which are monogenic in $\mathbb R^{m+1}$. 


\section{Construction of the sequence $\{\mathsf{M}_n^k(x)\}_{n\ge0}$}

Before starting with the construction of the sequence $\{\mathsf{M}_n^k(x)\}_{n\ge0}$, we shall first introduce some essential tools.

For a differentiable $\mathbb R$-valued function $\phi$ and a differentiable $\mathbb R_{0,m}$-valued function $g$, we have
\begin{equation}\label{Lr1}
\partial_{\underline x}(\phi g)=\partial_{\underline x}(\phi)g + \phi(\partial_{\underline x}g).
\end{equation}
Moreover, for a differentiable vector-valued function $\underline f=\sum_{j=1}^mf_je_j$, we also have
\begin{equation}\label{Lr2}
\partial_{\underline x}(\underline fg)=(\partial_{\underline x}\underline f)g-\underline f(\partial_{\underline x}g)-2\sum_{j=1}^mf_j(\partial_{x_j}g).
\end{equation}
Throughout this paper we denote by $\mathbf{P}_k(\underline x)$ an $\mathbb R_{0,m}$-valued homogeneous polynomial of degree $k\in\mathbb N_0$ which moreover is monogenic in $\mathbb R^m$, i.e. 
\begin{alignat*}{2}
\mathbf{P}_k(\underline x)&\in\mathbb R_{0,m},\quad\partial_{\underline x}\mathbf{P}_k(\underline x)=0,&\quad\underline x&\in\mathbb R^m,\\
\mathbf{P}_k(t\underline x)&=t^k\mathbf{P}_k(\underline x),\quad\underline x\in\mathbb R^m,&\quad t&\in\mathbb R.
\end{alignat*}
Let  
\[\beta_k(n)=\left\{\begin{array}{ll}n,&\text{if}\;n\;\text{even}\\2k+m+n-1,&\text{if}\;n\;\text{odd}\end{array}\right.\]
for $n\ge1$ and put $\beta_k(0)=1$. Using the Leibniz rules (\ref{Lr1})-(\ref{Lr2}) as well as Euler's theorem for homogeneous functions, we can deduce the useful equality:
\begin{equation}\label{ident1}
\partial_{\underline x}\big(\underline x^n\mathbf{P}_k(\underline x)\big)=-\beta_k(n)\underline x^{n-1}\mathbf{P}_k(\underline x),\quad n\ge1.
\end{equation}
One basic result in Clifford analysis is the so-called Cauchy-Kovalevskaya extension technique (see \cite{BDS,DSS}), which we will make heavy use in our paper. 

\begin{thm}\label{CK}
Every $\mathbb R_{0,m}$-valued function $g(\underline x)$ analytic in $\mathbb R^m$ has a unique monogenic extension $\mathsf{CK}[g]$ to $\mathbb R^{m+1}$, which is given by
\begin{equation}\label{CKf}
\mathsf{CK}[g(\underline x)](x)=\sum_{j=0}^\infty\frac{(-x_0)^j}{j!}\,\partial_{\underline x}^jg(\underline x).
\end{equation}
\end{thm}
Observe that a monogenic function $f(x)$ can be reconstructed by knowing its restriction to $\mathbb R^m$ using previous formula, i.e. 
\[f(x)=\mathsf{CK}[f(x)\vert_{\underline x_0=0}](x).\] 
It is also worth noting that
\begin{equation}\label{DCK}
\frac{1}{2}\,\overline\partial_x\mathsf{CK}[g(\underline x)](x)=-\partial_{\underline x}\mathsf{CK}[g(\underline x)](x)=\mathsf{CK}[-\partial_{\underline x}g(\underline x)](x).
\end{equation}
We are now ready to construct our sequence $\{\mathsf{M}_n^k(x)\}_{n\ge0}$ of $\mathbb R_{0,m}$-valued monogenic polynomials in $\mathbb R^{m+1}$ which satisfies the Appell condition (\ref{Acond}) and whose first term is $\mathsf{M}_0^k(x)=\mathbf{P}_k(\underline x)$.

It is easy to check that we can put 
\[\mathsf{M}_1^k(x)=\left(x_0+\frac{\underline x}{2k+m}\right)\mathbf{P}_k(\underline x),\;\mathsf{M}_2^k(x)=\left(x_0^2+\frac{2x_0\underline x}{2k+m}+\frac{\underline x^2}{2k+m}\right)\mathbf{P}_k(\underline x)\]
as the next two elements in our sequence. Indeed, they are monogenic in $\mathbb R^{m+1}$ and satisfy (\ref{Acond}). Thus, it seems that we can conjecture that each term $\mathsf{M}_n^k(x)$ in our sequence $\{\mathsf{M}_n^k(x)\}_{n\ge0}$ will be of the form
\[\mathsf{M}_n^k(x)=H_n(x_0,\underline x)\mathbf{P}_k(\underline x),\quad n\ge0,\]
where $H_n(x_0,\underline x)$ is a homogeneous polynomial of degree $n$ with real coefficients in the two variables $x_0$ and $\underline x$. Therefore
\[\mathsf{M}_n^k(x)\vert_{\underline x_0=0}=c_n\underline x^n\mathbf{P}_k(\underline x),\quad n\ge0,\]
for some real constant $c_n$ ($n\ge0$) with $c_0=1$. By Theorem \ref{CK}, we have that
\[\mathsf{M}_n^k(x)=\mathsf{CK}[\mathsf{M}_n^k(x)\vert_{\underline x_0=0}](x)=c_n\mathsf{CK}[\underline x^n\mathbf{P}_k(\underline x)](x),\quad n\ge0.\]
Since $\{\mathsf{M}_n^k(x)\}_{n\ge0}$ satisfy (\ref{Acond}), it follows from (\ref{DCK}) that
\[c_n\mathsf{CK}\big[-\partial_{\underline x}\big(\underline x^n\mathbf{P}_k(\underline x)\big)\big](x)=nc_{n-1}\mathsf{CK}[\underline x^{n-1}\mathbf{P}_k(\underline x)](x),\quad n\ge1.\]
Then, identity (\ref{ident1}) implies that 
\[c_n=\frac{nc_{n-1}}{\beta_k(n)},\quad n\ge1.\]
This recurrence relation is easy to solve. Indeed,
\[c_n=\frac{n!}{\displaystyle{\prod_{s=0}^n\beta_k(s)}},\quad n\ge0.\]
We thus get
\begin{equation}\label{fFt}
\mathsf{M}_n^k(x)=\frac{n!}{\displaystyle{\prod_{s=0}^n\beta_k(s)}}\,\mathsf{CK}[\underline x^n\mathbf{P}_k(\underline x)](x),\quad n\ge0.
\end{equation}
Moreover, using formula (\ref{CKf}) and identity (\ref{ident1}), we obtain
\begin{equation}\label{Mn(x)}
\mathsf{M}_n^k(x)=\left(\sum_{j=0}^n\binom{n}{j}C_{k,n}(j)x_0^j\underline x^{n-j}\right)\mathbf{P}_k(\underline x),\quad n\ge0,
\end{equation}
with $C_{k,n}(j)=\frac{\displaystyle{(n-j)!}}{\displaystyle{\prod_{s=0}^{n-j}\beta_k(s)}}$. 

\begin{prop}
Suppose that $\mathsf{M}_0^k(x)=\mathbf{P}_k(\underline x)$ is an $\mathbb R_{0,m}$-valued homogeneous polynomial of degree $k$ which is monogenic in $\mathbb R^m$. Then $\{\mathsf{M}_n^k(x)\}_{n\ge0}$ given by (\ref{Mn(x)}) is a sequence of $\mathbb R_{0,m}$-valued polynomials which are monogenic in $\mathbb R^{m+1}$ and satisfies the Appell condition (\ref{Acond}).
\end{prop}
\textbf{Remark:} \textit{For the particular case $k=0$, $\mathbf{P}_k(\underline x)=1$, $\{\mathsf{M}_n^k(x)\}_{n\ge0}$ equals the Appell sequence $\{\mathsf{P}_n^m(x)\}_{n\ge0}$ given by (\ref{ASMal}).}\vspace{0.15cm}

It should be noticed that $\{\mathsf{M}_n^k(x)\}_{n\ge0}$ is not a polynomial sequence. Indeed, $\mathsf{M}_n^k(x)$ ($n\ge0$) is a $\mathbb R_{0,m}$-valued homogeneous polynomial of degree $k+n$ which is monogenic in $\mathbb R^{m+1}$. 

It is worth pointing out that $\mathsf{M}_n^k(x)$ can be written as
\begin{equation}\label{axialmappell}
\mathsf{M}_n^k(x)=\Bigl(A_n(x_0,r)+\frac{\underline x}{r}\,B_n(x_0,r)\Bigr)\mathbf{P}_k(\underline x),\quad r=\vert\underline x\vert,\quad n\ge0,
\end{equation}
where $A_n$ and $B_n$ are $\mathbb R$-valued continuously differentiable functions in $\mathbb R^2$. The monogenicity of $\mathsf{M}_n^k(x)$ implies that $A_n$ and $B_n$ satisfy the Vekua-type system (see \cite{Ve}) 
\begin{equation*}
\left\{\begin{array}{ll}\partial_{x_0}A_n-\partial_rB_n&=\displaystyle{\frac{2k+m-1}{r}}\,B_n\\\partial_{x_0}B_n+\partial_rA_n&=0.\end{array}\right.
\end{equation*}
Monogenic functions of the form (\ref{axialmappell}) are called axial monogenic of degree $k$ and are an important class of functions in Clifford analysis (see \cite{LB,S1,S2}). Thus, we can also assert that $\{\mathsf{M}_n^k(x)\}_{n\ge0}$ is a sequence of axial monogenic functions of degree $k$ which satisfies the Appell condition (\ref{Acond}).


\section{Connection with Fueter's theorem}

Fueter's theorem constitutes a method for generating axial monogenic functions starting from holomorphic functions in the upper half of the complex plane. First proved by R. Fueter (see \cite{F}) in the setting of quaternionic analysis, this result was later extended to Clifford analysis in \cite{Sce,Q,S3}. For other references on this subject we refer the reader to e.g. \cite{KQS,LaRa,DS1,DS2,DS3,QS,Sp}.

For odd dimensions $m$, Fueter's theorem runs as follows (see e.g. \cite{S3}).

\begin{thm}
Let $f(z)=u(x,y)+iv(x,y)$ ($z=x+iy$) be a holomorphic function in some open subset $\Xi\subset\{z\in\mathbb C:\;y>0\}$. Put $\underline\omega=\underline x/r$, with $r=\vert\underline x\vert$. If $m$ is odd, then the function
\begin{equation*}
\mathsf{Ft}\left[f(z),\mathbf{P}_k(\underline x)\right](x)=\Delta_x^{k+\frac{m-1}{2}}\bigl[\bigl(u(x_0,r)+\underline\omega\,v(x_0,r)\bigr)\mathbf{P}_k(\underline x)\bigr]
\end{equation*}
is monogenic in $\widetilde\Omega=\{x\in\mathbb R^{m+1}:\;(x_0,r)\in\Xi\}$.
\end{thm}

In \cite{NGu} it was shown that for $m$ being odd the Appell sequence $\{\mathsf{P}_n^m(x)\}_{n\ge0}$ given by (\ref{ASMal}) and Fueter's theorem applied to the complex monomials $z^n$ are closely related. More precisely, $\mathsf{P}_n^m(x)$ equals (up to a multiplicative constant) $\mathsf{Ft}\left[z^{n+m-1},1\right](x)$ for $n\ge0$. Since for $k=0$ and $\mathbf{P}_k(\underline x)=1$, $\{\mathsf{M}_n^k(x)\}_{n\ge0}$ coincides with $\{\mathsf{P}_n^m(x)\}_{n\ge0}$, we may expect that a similar connection should exist between our sequence $\{\mathsf{M}_n^k(x)\}_{n\ge0}$ and Fueter's theorem. In fact, it was computed in \cite{DS1} (among other examples) that for odd dimensions $m$
\begin{multline}\label{ejFD}
\mathsf{Ft}\left[z^n,\mathbf{P}_k(\underline x)\right](x)=(-1)^{k+\frac{m-1}{2}}(2k+m-1)!!\alpha_k(n)\\
\times\mathsf{CK}[\underline x^{n-(2k+m-1)}\mathbf{P}_k(\underline x)](x),\quad n\ge2k+m-1,
\end{multline}
with
\[\alpha_k(n)=\left\{\begin{array}{ll}\displaystyle{\frac{n!!}{(n-2k-m+1))!!}},&\text{if}\;n\;\text{even}\\\alpha_k(n-1),&\text{if}\;n\;\text{odd},\end{array}\right.\]
and where $(\cdot)!!$ denotes the double factorial. We also note that 
\[\mathsf{Ft}\left[z^n,\mathbf{P}_k(\underline x)\right](x)=0,\quad\text{for}\;n<2k+m-1.\]
Thus, using (\ref{fFt}) and (\ref{ejFD}), we get:

\begin{prop}
For odd dimensions $m$, each term $\mathsf{M}_n^k(x)$ of the sequence $\{\mathsf{M}_n^k(x)\}_{n\ge0}$ equals (up to a multiplicative constant) $\mathsf{Ft}\left[z^{n+2k+m-1},\mathbf{P}_k(\underline x)\right](x)$.
\end{prop}


\subsection*{Acknowledgment}

Financial support from \emph{Funda\c{c}\~ao para a Ci\^encia e a Tecnologia} (FCT), Portugal (Post-Doctoral Grant Number: SFRH/BPD/45260/2008), is gratefully acknowledged. The author is grateful to F. Sommen for teaching him those techniques which are essential to obtain the results presented in this paper. He also wishes to thank I. Ca\c{c}\~ao and H. R. Malonek for useful conversations on the topic of Appell sequences.    


\end{document}